\documentclass[11pt, a4paper, twoside]{article}
\usepackage[a4paper]{anysize}
\usepackage{amsmath, amsfonts, amssymb, amsthm, graphicx}
\usepackage{amsmath,amssymb}
\usepackage[latin1]{inputenc}

\usepackage{amssymb}
\usepackage{amsmath}
\usepackage{graphics}
\usepackage{latexsym}
\usepackage{amsfonts}
\usepackage{color}
\usepackage{epstopdf}
\usepackage{ctable}
\usepackage{graphicx}
\usepackage{verbatim}
\newtheorem{theo}{Theorem}[section]
\newtheorem{lemma}[theo]{Lemma}

\newtheorem{cor}[theo]{Corollary}
\newtheorem{prop}[theo]{Proposition}
\newtheorem{uppg}[theo]{Exercise}

\newcommand{\be}{\begin{eqnarray*}}
\newcommand{\ee}{\end{eqnarray*}}
\newcommand{\ben}{\begin{eqnarray}}
\newcommand{\een}{\end{eqnarray}}

\def\subsecn (#1) {\medskip\ \ \ {\it #1}\medskip}

\newcommand{\lp}[1]{\left(\begin{array}{#1}}
\newcommand{\rp}{\end{array}\right)}

\newcommand{\leftd}[1]{\left\{\begin{array}{#1}}
\newcommand{\rightd}{\end{array}\right.}

\def\A {\mathbf{A}}
\def\B {\mathbf{B}}

\def\D {\mathbf{D}}
\def\E {\mathbf{E}}

\def\I {\mathbf{I}}

\def\P {\mathbf{P}}

\def\R {\mathbf{R}}

\def\a {\boldsymbol{a}}
\def\b {\boldsymbol{b}}

\def\e {\boldsymbol{e}}

\def\u {\boldsymbol{u}}
\def\v {\boldsymbol{v}}

\def\x {\boldsymbol{x}}
\def\y {\boldsymbol{y}}

\def\Eb {\mathbb{E}}

\date{\today}
\begin{document}
\title{Oscillation of adaptative Metropolis-Hasting and simulated annealing algorithms around  penalized least squares estimator\\
\small Azzouz Dermoune, Daoud Ounaissi, Nadji Rahmania\\
Laboratoire Paul Painlev\'e, USTL-UMR-CNRS 8524.\\
UFR de Math\'ematiques, B\^at. M2, 59655 Villeneuve d'Ascq C\'edex, France.\\
azzouz.dermoune@univ-lille1.fr}

\maketitle
{\bf Abstract.} 
In this work we study, as the temperature goes to zero, the oscillation of Metropolis-Hasting's algorithm around the Basis Pursuit De-noising solutions.  
We derive new criteria for choosing the proposal distribution and the temperature in Metropolis-Hasting's algorithm. Finally we apply these results to compare Metropolis-Hasting's 
and  simulated annealing algorithms.

{\bf keyword.} Penalized least squares. Adaptative Metropolis-Hasting. Simulated annealing algorithms. Gibbs measures.
 
\section{Penalized least squares estimate} 
Let $\A$ and $\y$ be respectively an $n\times p$ measurement matrix and a $n\times 1$ measurement vector. 
The unknown vector $\x$ belongs to $\R^p$. We are interested in the case where the number of parameters $p$ is larger than the data number $n$. Given the penalty function $\x\to \|\x\|_1:=\sum_{i=1}^p|x_i|$ and the smoothing parameter $t\geq 0$, the penalized least squares estimate (PLSE in short)
proposes to recover the vector $\x$ using the minimization problem   
$\x(\y,t)\in\arg\min\{\|\x\|_1+\frac{\|\A\x-\y\|^2}{2t}:\quad \x\in\R^p\}$ (known as Basis Pursuit De-Noising method \cite{Donoho}).
Here 
$\|\cdot\|$ denotes the Euclidean norm. The set of PLSE can be found using (FISTA) algorithm \cite{Beck}.
In our work we consider the family of probabilities (called also Gibbs measures)
\ben 
P_T^{\y,t}:=\frac{\exp(-\frac{1}{T}F(\x,\y,t))d\x}{\int_{\R^p}\exp(-\frac{1}{T}F(\x,\y,t))d\x},
\label{PT}
\een 
where $T >0$ is called the temperature and $F(\x,\y,t)=\|\x\|_1+\frac{\|\A\x-\y\|^2}{2t}$ is called the objective function.      
Well known results tell us that the family of the probabilities  (\ref{PT}) 
oscillates around the set of PLSE as $T\to 0$.  
More precisely, any sequence $(\P_{T_k}^{\y,t}: T_k\to 0)$ is tight \cite{Hwang}, \cite{AthreyaHwang} i.e. we can extract a convergent subsequence from 
$(P_{T_k}^{\y,t})$. If $P_{T_k}^{\y,t}\to P^{\y,t}$, then $P^{\y,t}$ concentrates on $\arg\min\{F(\x,\y,t):\quad \x\in\R^p\}$. Hence, using Metropolis-Hasting's algorithm with small temperature and the target (\ref{PT}) or the  
simulated annealing algorithm, we can construct Markov chains having the tails
located near the set of PLSE. Fort et al., in    
a recent work \cite{Moulines}, propose a new algorithm based on Metropolis and Langevin 
equation.       

The efficiency of Metropolis-Hasting and simulated annealing algorithms depends on the choice of the proposal distribution and 
the temperature. In Section 2 we give a precise scaling of the asymptotic of the measures (\ref{PT}) 
as $T\to 0$. In Section 3 we derive new criteria of the choice of the proposal distribution and the temperature. We also apply  these  criteria to compare  Metropolis-Hasting  
and the simulated annealing algorithms. Finally we numerically illustrate our results .  
  
\section{Gibbs measures scaling as the temperature goes to zero}  

First, we need some notations
The vector $sgn(\x)$ will denotes the $p$ by 1 matrix with the components 
$sgn(x_i)=1$ if $x_i >0$, $sgn(x_i)=-1$ if $x_i <0$ and $sgn(0)$ is any element of $[-1,1]$. 
We will denote, for each subset $I\subset\{1, \ldots, p\}$ and for each vector $\v\in \R^p$, $\v(I)=(v(i): i\in I)\in \R^I$. The notation $v\leq w$ means $v(i)\leq w(i)$ for all $i$. The scalar product is denoted by $\langle\cdot, \cdot\rangle$, 
and $(\e_i:\quad i=1, \ldots)$ denotes the canonical basis of $\R^p$.  

Let us recall some properties of the Basis Pursuit De-noising   
minimizers.
\begin{prop}\label{PR} A vector $\x(\y,t)$ is a minimizer of  the map 
$\x\to \|\x\|_1+\frac{\|\A\x-\y\|^2}{2t}$ if the vector $\xi(\y,t)=\frac{\A^*(\y-A\x(\y,t))}{t}$
belongs to $sgn(\x(\y,t))$. The vectors  
$\xi(\y,t)$, $\A\x(\y,t)$ and the l1-norm $\|\x(\y,t)\|_1$ are constant on the set of PLSE.
Here $\A^*$ denotes the transpose of the matrix $\A$. 

\end{prop}  
The sets    
$I_0=\{i\in\{1,\ldots, p\}:\quad \x_i(\y,t)=0\}$, $\partial I_0=\{i\in I_0:\quad |\xi_i(\y,t)|=1\}$ will play an important role in the Gibbs measures scaling.    
The set $S=\{1,\ldots, p\}\setminus I_0$ is the
support of the PLSE $\x(\y,t)$ i.e. $S=\{i\in \{1, \ldots, p\}: \x_i(\y,t)\neq 0\}$.
In the sequel $X_T(\y,t)$ will denote a random vector having the probability distribution (\ref{PT}). If the set of PLSE is a singleton $\x(\y,t)$, then we can show that 
$X_T(\y,t)\to \x(\y,t)$ in probability as $T\to 0$ see e.g. \cite{AthreyaHwang}.  

Before announcing our main result we need some preliminary lemmas. 
\begin{lemma}\label{Fsimplification} 
Let $\x(\y,t)$ be any PLSE and $m(\y,t)=F(\x(\y,t),\y,t)$ be the minimum of the objective function 
$F(\x,\y,t)$. The function $F(\x,\y,t)-m(\y,t)$ is equal to  
\ben\label{Fx} 
\sum_{i=1}^p|x_i|(1-sgn(x_i)\xi_i(\y,t))+\frac{\|\A(\x-\x(\y,t))\|^2}{2t}.
\een 
If $\x$ is near the PLSE $\x(\y,t)$, then $F(\x,\y,t)-m(\y,t)$ becomes  
\ben\label{Fxnear} 
\sum_{i\in I_0}|x_i|(1-sgn(x_i)\xi_i(\y,t))+\frac{\|\A(\x-\x(\y,t))\|^2}{2t}.  
\een 
\end{lemma} 
{\bf Proof.} From the equality $\|\A\x-\y\|^2=\|\A(\x-\x(\y,t))\|^2+2\langle\A(\x-\x(\y,t)),\A\x(\y,t)-\y\rangle+
\|\A\x(\y,t)-\y\|^2$, we have  
\be
&&F(\x,\y,t)=\\ 
&&\|\x\|_1+\frac{\|\A(\x-\x(\y,t))\|^2}{2t}+\frac{\langle\A(\x-\x(\y,t)),\A\x(\y,t)-\y\rangle}{t}+\frac{\|\A\x(\y,t)-\y\|^2}{2t}\\
&&=\|\x\|_1+\frac{\|\A(\x-\x(\y,t))\|^2}{2t}+\frac{\langle\x-\x(\y,t),\A^*(\A\x(\y,t)-\y)\rangle}{t}+\frac{\|\A\x(\y,t)-\y\|^2}{2t}.
\ee 
From the equation 
$\xi(\y,t)=\frac{\A^*(y-\A\x(\y,t))}{t}$ Proposition (\ref{PR}), we have 
\ben\label{fc} 
\frac{\langle\x-\x(\y,t),\A^*(\A\x(\y,t)-\y)\rangle}{t}&=&-\langle\x-\x(\y,t),\xi(\y,t)\rangle\nonumber\\ 
&=&-\langle \x,\xi(\y,t)\rangle+\|\x(\y,t)\|_1.
\een 
Now formulas (\ref{Fx}) and (\ref{Fxnear}) are an easy consequence of the formula (\ref{fc}).

The following lemma gives a sufficient condition for the uniqueness of 
the PLSE $\x(\y,t)$. 
\begin{prop}\label{uniqueness} 
If the matrix $[\langle \A\e_i,\A\e_j\rangle,\quad i,j\in (S\cup\partial I_0)]$ 
is invertible, then the set of PLSE is a singleton. 
\end{prop}  
{\bf Proof.} Observe that the invertibility of the matrix $[\langle \A\e_i,\A\e_j\rangle,\quad i,j\in (S\cup\partial I_0)]$ 
is equivalent to say that the linear operator 
$\A_{S\cup\partial I_0}:\R^{S\cup\partial I_0}\to \R^n$ is injective. Here $\A_{S\cup\partial I_0}$ denotes the sub-matrix
of $\A$ having the columns indexed by $S\cup\partial I_0$. The inverse of $\A_{S\cup\partial I_0}$ defined from $\R^{S\cup\partial I_0}$ into its range  $R(\A_{S\cup\partial I_0})$ is denoted by $\A_{S\cup\partial I_0}^{-1}$.   
Now, we recall a result of Grasmair et al. \cite{Grasmair} Lemma 3.10. 
Let $M(\x(\y,t)):=\max\{|\xi_i(\y,t))|:\quad i\in I_0\setminus \partial I_0\}$, and for any couple $\x^{(1)}, \x^{(2)}\in\R^p$, 
$D(\x^{(1)},\x^{(2)}):=\|\x^{(1)}\|_1-\|\x^{(2)}\|_1-\langle\xi, \x^{(1)}-\x^{(2)}\rangle$
for some fixed $\xi\in sgn(\x^{(2)})$. 
The result of Grasmair et al. 
tells us that for all $x$,  
\be 
\|\x-\x(\y,t)\|\leq \|\A_{S\cup\partial I_0}^{-1}\|\|\A (\x-\x(\y,t))\|+\frac{1+\|\A_{S\cup\partial I_0}^{-1}\|\|\A\|}{1-M(\x(\y,t))}
D(\x,\x(\y,t)), 
\ee 
where $\|\B\|$ denotes the operator norm of the matrix $\B$. If $\x$ is another PLSE, then from Proposition (\ref{PR}), we have 
$\A\x=\A\x(\y,t)$ and  
$D(\x,\x(y,t))=0$, which achieves the proof.

Now we can announce our last lemma. 
 \begin{lemma}\label{evenementnegligeable} Let $\partial I_0=K_1\cup K_2$ be a partition such that $K_1, K_2\neq \emptyset$ and 
$E(K_1,K_2)=E_{-1}(K_1)\cap E_1(K_2)$ with 
$E_{-1}(K_1)=\{\x\in\R^p:\, sgn(x_i)\xi_i(\y,t)\\=-1,\,\forall\,i\in K_1\}$, and 
$E_{+1}(K_2)=\{\x\in\R^p:\, sgn(x_i)\xi_i(\y,t)=1,\,\forall\,i\in K_2\}$. 
If $[\langle \A\e_i,\A\e_j\rangle,\quad i,j\in (S\cup\partial I_0)]$ 
is invertible, then the set of PLSE is a singleton and   
the probability of the event $E_T(K_1,K_2):=[X_T(\y,t)\in E(K_1,K_2)]$ 
tends to 0 as $T\to 0$. As a consequence, we have 
$\P(E_T(\emptyset,\partial I_0))\to 1$ as $T\to 0$.  
\end{lemma}  
{\bf Proof.} The uniqueness of the PLSE is shown in the Proposition (\ref{uniqueness}). 
Now, we prove the rest of our Lemma.   
We have $\P(X_T(\y,t)\in E(K_1,K_2))=\frac{A_T(K_1,K_2)}{B_T}$, 
where $A_T(K_1,K_2)=\int_{E(K_1,K_2)}\exp(-\frac{1}{T}F(\x,\y,t))d\x$, and 
$\int\exp(-\frac{1}{T}F(\x,\y,t))d\x\\=B_T$.
We know that for small $T$, $X_T(\y,t)$ will concentrate on $\x(\y,t)$.
It follows that the PDF (\ref{PT}) becomes more and more concentrated near $\x(\y,t)$. Hence, it is sufficient to consider, for small $\delta$,   
\be 
A_T(K_1,K_2,\delta)=\int_{E(K_1,K_2,\delta)}\exp(-\frac{1}{T}F(\x,\y,t))d\x,\\
B_T(\delta)=\int_{\|\x-\x(\y,t)\|_{\infty}\leq \delta}\exp(-\frac{1}{T}F(\x,\y,t))d\x, 
\ee
where $E(K_1,K_2,\delta)=E(K_1,K_2)\cap [\x: \|\x-\x(\y,t)\|_{\infty}\leq \delta]$ and 
$\|\x\|_{\infty}=\max(|x_i|: i=1, \ldots, p)$. 
From the Lemma (\ref{Fsimplification}) formula (\ref{Fxnear}), we have 
\be  
&&A_T(K_1,K_2,\delta)=\exp(-\frac{m(\y,t)}{T})\\
&&\int_{E(K_1,K_2,\delta)}\exp(-\frac{1}{T}(\sum_{l\in I_0}|x_l|(1-sgn(x_l)\xi_l(\y,t))+\frac{\|\A(\x-\x(\y,t))\|^2}{2t}) )d\x\\
&&=\exp(-\frac{m(\y,t)}{T})\int_{E(K_1,K_2,\delta)}\\
&&\exp(-\frac{1}{T}(\sum_{l\in I_0\setminus K_2}|x_l|(1-sgn(x_l)\xi_l(\y,t))+\frac{\|\A(\x-\x(\y,t))\|^2}{2t})d\x.
\ee 
Using the change of variables 
\ben 
\u=\frac{\x(I_0\setminus K_2)}{T},\quad \frac{\x(S\cup K_2)-\x(\y,t,S\cup K_2)}{\sqrt{T}}=\v,
\label{CV}
\een 
we get $A_T(K_1,K_2,\delta)=\exp(-\frac{m(\y,t)}{T})T^{\frac{|I_0\setminus K_2|+p}{2}}
C_T(K_1,K_2)$, where 
$C_T(K_1,K_2)=\int_{\tilde{E}_T(K_1,K_2,\delta)}\\
\exp(-(\sum_{i\in I_0\setminus K_2}|u_i|(1-sgn(u_i)\xi_i(\y,t))+\frac{\|\sqrt{T}
\sum_{i\in I_0\setminus K_2}u_i\A\e_i+\sum_{i\in (S\cup K_2)}v_i\A\e_i\|^2}{2t})d\u d\v$,  
and 
\be 
&&\tilde{E}_T(K_1,K_2,\delta)=\{\u\in [-\frac{\delta}{T},\frac{\delta}{T}]^{I_0\setminus K_2},\v\in [-\frac{\delta}{\sqrt{T}},\frac{\delta}{\sqrt{T}}]^{S\cup K_2}:\\
&&sgn(\u_{K_1})=-\xi_{K_1}(\y,t),\,
sgn(\v_{K_2})=\xi_{K_2}(\y,t),\\
&& sgn(\sqrt{T}\v_S+\x_{S}(\y,t))=sgn(\x_{S}(\y,t))\}
\ee
and $|I|$ denotes the cardinality of the set $I$. 
From the same calculation we can show that 
$B_T(\delta)=\sum_{K_1',K_2':\partial I_0=K_1'\cup K_2'}A_T(K_1',K_2',\delta)$.   
We emphasize that the couple $K_1'=\emptyset, K_2'=\partial I_0$ is an element of the latter sum. 
Moreover, the quantity $\frac{p+|I_0\setminus K_2|}{2}$
is minimal at $K_2=\partial I_0$. From this  we derive that 
\be 
\frac{A_T(K_1,K_2,\delta)}{B_T(\delta)}=
\frac{T^{\frac{|I_0\setminus K_2|-|I_0\setminus\partial I_0|}{2}}C_T(K_1,K_2)}
{C_T(\emptyset,\partial I_0)+
\sum_{K_1',K_2'\neq \emptyset :\partial I_0=K_1'\cup K_2'}
T^{\frac{|I_0\setminus K_2'|-|I_0\setminus\partial I_0|}{2}}C_T(K_1',K_2')}
\ee 
converges to 0 as $T\to 0$, because $C_T(K_1',K_2')\to C_0(K_1',K_2')\neq 0$ as $T\to 0$
for any  partition $K_1', K_2'$ of $\partial I_0$.

Our new criteria of the choice of the proposal distribution and the temperature in Metropolis-Hasting and the simulated annealing algorithms are based on the following result. 
\begin{prop} \label{mainpropo} 
Suppose that the matrix 
$[\langle \A\e_i,\A\e_j\rangle,\quad i,j\in (S\cup \partial I_0)]$ 
is invertible. Then the random vector 
$(\frac{X_T(\y,t,i)}{T}:\quad i\in (I_0\setminus\partial I_0)), 
(\frac{X_T(\y,t,i)-x(\y,t,i)}{\sqrt{T}}:\quad i\in (S\cup \partial I_0))$
converges to the random vector 
$(X_i(\y,t):\quad i\in (I_0\setminus\partial I_0)), (X_i(\y,t):\quad i\in (S\cup \partial I_0))$
having the PDF proportional to  
\be 
&&\prod_{i\in (I_0\setminus \partial I_0)} \exp(-|x_i|(1-sgn(x_i)\xi_i(\y,t))\\
&&\exp(-\frac{\|\sum_{i\in (S\cup\partial I_0)}x_i\A\e_i\|^2}{2t})\prod_{i\in\partial I_0}1_{[sgn(x_i)\xi_i(\y,t))=1]}. 
\ee 
\end{prop} 
{\bf Proof.} 
Let $I=I_0\setminus \partial I_0$ and $J=S\cup \partial I_0$ and $\a, \b\in\R^p$. We want to prove that  
$\P(\a(I)\leq \frac{X_T(\y,t,I)}{T}\leq \b(I),\a(J)\leq \frac{X_T(\y,t,J)-\x(\y,t,J)}{\sqrt{T}}
\leq \b(J))$ converges to 
$\P(\a(I)\leq X(\y,t,I)\leq \b(I),\a(J)\leq X(\y,t,J)-\x(\y,t,J)
\leq \b(J))$ as $T\to 0$. 
As we shown in the Lemma (\ref{evenementnegligeable}), it is sufficient to consider, for small $\delta$,  
\be 
&&\P(\a(I)\leq \frac{X_T(\y,t,I)}{T}\leq \b(I),\a(J)\leq \frac{X_T(\y,t,J)-\x(\y,t,J)}{\sqrt{T}}
\leq \b(J),\\
&&\|X_T(\y,t)-\x(\y,t)\|_{\infty} \leq \delta)\\
&&=\sum_{K_1,K_2:\partial I_0=K_1\cup K_2}\P(\ldots\,|\,E_T(K_1,K_2,\delta))\P(E_T(K_1,K_2,\delta)),
\ee 
where the events $E_T(K_1,K_2)$ are defined in the Lemma (\ref{evenementnegligeable}). As we are interested in the limit 
as $T\to 0$ and thanks to the lemma (\ref{evenementnegligeable}) only the term 
$P(\cdots\,|\,E_T(\emptyset,\partial I_0,\delta))\P(E_T(\emptyset,\partial I_0,\delta))$ is needed.  
More precisely we have only to study the term  
\be 
&&\P(\a(I)\leq \frac{X_T(\y,t,I)}{T}\leq \b(I),\\
&&\a(J)\leq \frac{X_T(\y,t,J)-\x(\y,t,J)}{\sqrt{T}}
\leq \b(J)\,|\,E_T(\emptyset,\partial I_0,\delta)
=\frac{A_T(\delta)}{B_T(\delta)}
\ee
where 
\be 
&&A_T(\delta)=\int_{T\a(I)}^{T\b(I)}\int_{\sqrt{T}\a(J)+\x(\y,t,J)}^{\sqrt{T}\b(J)+\x(\y,t,J)}\exp(-\frac{1}{T}F(\x,\y,t))1_{E(\emptyset,\partial I_0,\delta)}(\x)d\x\\
&&B_T(\delta)=\int\exp(-\frac{1}{T}F(\x,\y,t))1_{E(\emptyset,\partial I_0,\delta)}d\x.
\ee 
From the Lemma (\ref{Fsimplification}) we have 
\be  
&&A_T(\delta)=\exp(-\frac{m(\y,t)}{T})\int_{T\a(I)}^{T\b(I)}\int_{\sqrt{T}\a(J)+\x(\y,t,J)}^{\sqrt{T}\b(J)+\x(\y,t,J)}\\
&&\exp(-\frac{1}{T}(\sum_{i\in I}|x_i|(1-sgn(x_i)\xi_i(\y,t))+\frac{\|\A(\x-\x(\y,t))\|^2}{2t}) )1_{E(\emptyset,\partial I_0,\delta)}(\x)d\x
\ee 
Using the change of variables 
\ben 
\u=\frac{\x(I)}{T},\quad \frac{\x(J)-\x(\y,t,J)}{\sqrt{T}}=\v,
\label{CV}
\een 
we get 
\be 
&&A_T(\delta)=\exp(-\frac{m(\y,t)}{T})T^{\frac{|I|+p}{2}}\int_{\a(I)}^{\b(I)}\int_{\a(J)}^{\b(J)}d\u d\v\\
&&\exp(-(\sum_{i\in I}|u_i|(1-sgn(u_i)\xi_i(\y,t))+\frac{\|\sqrt{T}\sum_{i\in I}u_i\A\e_i+\sum_{i\in J}v_i\A\e_i\|^2}{2t}))\\
&&1_{\tilde{E}_T(\emptyset,\partial I_0,\delta)}(\u,\v). 
\ee    
Now we are going to study $T^{-\frac{|I|+p}{2}}B_T\exp(\frac{m(\y,t)}{T})$. From  the change of variables formula (\ref{CV}), 
we have 
\be 
&&\lim_{T\to 0} T^{-\frac{|I|+p}{2}}B_T(\delta)\exp(\frac{m(\y,t)}{T})=\lim_{T\to 0}T^{-\frac{|I|+p}{2}}\int_{-\delta\leq \x\leq \delta}\\
&&\exp(-\frac{1}{T}(\sum_{i\in I}|x_i|(1-sgn(x_i)\xi_i(\y,t))+\frac{\|\A(\x-\x(\y,t))\|^2}{2t}) )1_{E(\emptyset,\partial I_0,\delta)}d\x\\
&&=\int\exp(-(\sum_{i\in I}|u_i|(1-sgn(u_i)\xi_i(\y,t))+\frac{\|\sum_{i\in J}v_i\A\e_i\|^2}{2t}))\\
&&\prod_{i\in\partial I_0}1_{[sgn(v_i)\xi_i(\y,t)=1]}dudv,
\ee  
which achieves the proof.

\section{One dimensional case} 
In the one dimensional case the objective function 
$F(x,y,t)=|x|+\frac{(x-y)^2}{2t}$. 
In this case 
$x(y,t)=0$, for $|y| \leq t$, $ x(y,t)=y+t$, for $y < -t$, and  
$x(y,t)=y-t$, for $y > t$.  

Let $X_{T}(y,t)$ be a random variable drawn from the PDF proportional to 
$\exp\left(-\frac{1}{T}(|x|+\frac{(x-y)^2}{2t})\right)$.
The following is a consequence of Proposition (\ref{mainpropo}) and precise, for $y >0$, the behavior of $X_{T}(y,t)$. 
\begin{prop}\label{toyresults} 1) If $y\in [0,t)$, then $\frac{X_{T}(y,t)}{T}\to X(y,t)$, where 
$X(y,t)$ is the random variable having the PDF 
\be 
x\to \frac{1-\frac{y^2}{t^2}}{2}\exp(-|x|(1-sgn(x)\frac{y}{t})).
\ee 
2) Known the event $[X_T(t,t)< 0]$, the random variable $\frac{X_T(t,t)}{T}\to -\mathcal{E}(2)$
where $\mathcal{E}(2)$ is the random variable having the exponential distribution with the parameter 2, i.e. the PDF of $\mathcal{E}(2)$ is equal to 
$2\exp(-2x)1_{[x > 0]}$.\\ 
3) Known the event $[X_T(t,t) >0]$, the random variable $\frac{X_T(t,t)}{\sqrt{T}}\to |N(0,t)|$, 
where $N(0,t)$ is the standard Gaussian with the variance $t$.\\ 
4) We have for $y >t$ that 
$\frac{X_{T}(y,t)-(y-t)}{\sqrt{T}}\to N(0,t)$
as $T\to 0$.
\end{prop} 
The following corollary is a simple case of the lemma (\ref{evenementnegligeable}). 
\begin{cor} We have $\P(X_T(t,t)<0)\to 0$ as $T\to 0$. It follows that  
and $\frac{X_T(t,t)}{\sqrt{T}}$
converge to $|N(0,t)|$. Roughly speaking $X_T(t,t)\approx \sqrt{T}|N(0,t)|$ as
$T\to 0$.
\end{cor}

\subsection{Interpretation of Proposition (\ref{toyresults})} 
If $0\leq y < t$, then the density of $X(y,t)$ is a mixture of exponential probability distributions i.e. is equal to 
$\frac{1-\frac{y}{t}}{2}f_{X_-(y,t)}(x)+\frac{1+\frac{y}{t}}{2}f_{X_+(y,t)}(x)$,  
where $X_-(y,t), X_+(y,t)$ are independent variables having respectively the exponential distribution 
$-\mathcal{E}(1+\frac{y}{t}), \mathcal{E}(1-\frac{y}{t})$. Hence, $X(y,t)$ has the same PDF as
$X_{b(\frac{1+\frac{y}{t}}{2})}(y,t)$, where $(X_-(y,t),X_+(y,t),b(\frac{1+\frac{y}{t}}{2}))$ are independent with  the PDF  
\be 
-\mathcal{E}(1+\frac{y}{t}),&& \mathcal{E}(1-\frac{y}{t}),\\
\P(b(\frac{1+\frac{y}{t}}{2})=-)=\frac{1-\frac{y}{t}}{2}, &&\P(b(\frac{1+\frac{y}{t}}{2})=+)=\frac{1+\frac{y}{t}}{2}
\ee 
respectively. 
We know, for $y\in (0,t)$, that $X_T(y,t)$ converges to the Dirac measure $\delta_0$.
Hence, we have  for small $T$ that $X_T(y,t)\approx \delta_0$. 
Proposition (\ref{toyresults}) makes a zoom on the latter convergence. 
It shows for $y\in (0,t)$ and small $T$ that 
$X_T(y,t)\approx TX_{b(\frac{1+\frac{y}{t}}{2})}(y,t)$  
and shows that $X_T(y,t)-(y-t)\approx \sqrt{T}\mathcal{N}(0,t)$ for $y\geq t$. 

Using this approximation we will discuss how the proposal distribution 
in Metropolis-Hasting's depends on the data $y, t$ and the temperature $T$. 
We will also discuss the choice of the temperature in 
the simulated annealing algorithm. In 
Figure 1 we plot the probability density function of $X(y,t)$ when $y\in (0,t)$. 
\begin{figure}[!ht]
  \centering
     \includegraphics[width=8cm]{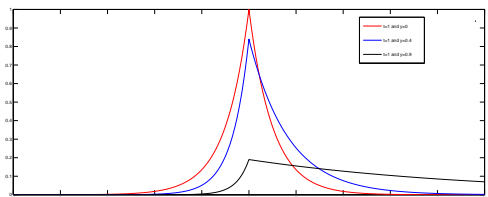}
 \caption{The density of $X(y,1)$ for $y=0, 0.4, 0.9$.}
\end{figure}\\

 \section{Numerical results} 
 \subsection{Choosing the proposal distribution in Metropolis-Hasting's algorithm} 
We want to sample from $X_T(y,t)$ using Metropolis-Hasting's algorithm with a family 
of proposal distributions. 
There are many criteria to choose the best proposal distribution 
see e.g. 
\cite{Johansen} example 5.3 chapter 5 and 
Gelman et al. \cite{Gelman}. In the sequel we propose new criteria based on the asymptotic distribution given in Proposition (\ref{toyresults}).   
We distinguish three cases. 

1) The case $y\in (0,t)$. 
 
a) Criterion using the asymptotic bias: We propose $X_T(y,t)$ as an estimator of $soft(y,t)=0$. Its bias, for small $T$, is equal to 
\be 
\Eb[X_T(y,t)]\approx T\Eb[X(y,t)]=T\frac{2\frac{y}{t}}{(1-\frac{y^2}{t^2})}
=Tm_1(\frac{y}{t}).
\ee  
The best proposal for sampling $X_T(y,t)$ will produce a sequence $(\theta^{(n)}(y):\quad n=1, \ldots, N)$ such that $\frac{1}{N}\sum_{n=1}^N\theta^{(n)}(y)$ 
is the nearest to $Tm_1(\frac{y}{t})$. In order to take account of all $y\in (0,t)$ we consider a sample 
$(U_i: i=1, \ldots M)$  of the $Beta(\alpha,\beta)$ distribution  with $\alpha=1$ and $\beta=3$ . For each proposal $q$ we calculate 
the objective function 
$f_1(q)=\frac{1}{M}\sum_{i=1}^M|\frac{1}{N}\sum_{n=1}^N\theta^{(n)}(tU_i)-Tm_1(U_i)|$.
We say that the proposal $q^*$ is the best among a family $F$ 
of proposal distributions if 
$q^*$ is the minimizer of $q\in F\to f_1(q)$. We tried others parameters of Beta distribution and also 
Uniform distribution on $(0,t)$. We showed that our criterion is unstable for these choices.       

b) Criterion using the asymptotic mean square error: The mean square error $\Eb[X_T^2(y,t)]$ for small $T$ is equal to 
\be 
T^2\Eb[X^2(y,t)]=T^2\{\frac{1-\frac{y}{t}}{(1+\frac{y}{t})^2}+\frac{1+\frac{y}{t}}{(1-\frac{y}{t})^2}\}=T^2m_2(\frac{y}{t}). 
\ee  
Now we can announce our second criterion. 
The best proposal for sampling $X_T(y,t)$ will produce a sequence $(\theta^{(n)}:\quad n=1, \ldots, N)$ such that 
$\frac{1}{N}\sum_{n=1}^N(\theta^{(n)})^2$
is the nearest to $T^2m_2(\frac{y}{t})$ for all $y\in (0,t)$. Similarly to a), we propose for any family $F$ of proposal distributions, the best proposal distribution as the minimizer of 
\be 
q\in F\to f_2(q)=\frac{1}{M}\sum_{i=1}^M|\frac{1}{N}\sum_{n=1}^N(\theta^{(n)}(tU_i))^2-T^2m_2(U_i)|.
\ee  
If $\arg\min_{F} f_1\neq \arg\min_{F}f_2$, then we propose the minimizer of 
$q\in F\to f_1(q)+f_2(q)$ as the best proposal distribution. \\
In order to illustrate these results we consider  $M=600$ chains with size $N=5000$,  with the proposal distribution $\mathcal{N}(0,\sigma^2)$ with different values of $\sigma^2$. The table 1 shows that the best proposal distribution for $t=1$ is $\mathcal{N}(0,1)$. 
\\
\begin{table}[!ht]
\begin{center}
  \begin{tabular}{|c||ccc|}
    \hline
    \scriptsize{Proposal} & \scriptsize{$f_1(q)$} &  \scriptsize{$f_2(q)$}  &  \scriptsize{$f_1(q)+f_2(q)$} \\
    \hline
    \scriptsize{$\sigma^2=1$} & \scriptsize{0.0351} & \scriptsize{0.0615} & \scriptsize{0.0966} \\ 
     \hline
    \scriptsize{$\sigma^2=9$ } & \scriptsize{0.0373}  & \scriptsize{0.0605} & \scriptsize{0.0978} \\
      \hline
    \scriptsize{$\sigma^2=16$ } & \scriptsize{0.0394} & \scriptsize{0.0604} & \scriptsize{0.0998} \\
      \hline
\end{tabular}
  \caption{ $y<t$, $t=1$, $T=0.1$,  $N=5000$, $M=600$ .}
  \end{center}
\end{table}\\
2) The case $y=t$. We showed that for small $T$ the random variable 
$X_T(t,t)$ is approximatly equal to $\sqrt{T}|\mathcal{N}(0,t)|$.   

a) Criterion using the asymptotic bias: The bias of $X_T(t,t)$, for small $T$, is equal to 
$\Eb[X_T(t,t)]\approx \sqrt{T}\Eb[|\mathcal{N}(0,t)|]=\sqrt{\frac{2Tt}{\pi}}$. 
The best proposal distribution $q$ for sampling positive values of $X_T(t,t)$ will produce a sequence $(\theta^{(n)}:\quad n=1, \ldots, N)$ such that 
\be 
f_1(q)=|\frac{1}{card\{n\leq N:\quad\theta^{(n)}> 0\}}\sum_{n=1}^N\theta^{(n)}1_{[\theta^{(n)}> 0]}-\sqrt{\frac{2Tt}{\pi}}|
\ee 
is minimal. 

b) Criterion using the asymptotic mean square: The mean square error $\Eb[X_T^2(t,t)]$ for small $T$ is equal to 
$T\Eb[\mathcal{N}^2(0,t)]=Tt$.  
The best proposal distribution $q$ for sampling $X_T(t,t)$ will produce a sequence $(\theta^{(n)}:\quad n=1, \ldots, N)$ such that 
\be 
f_2(q)=|\frac{1}{N}\sum_{n=1}^N(\theta^{(n)})^2-Tt|
\ee
is minimal. If the minimizers of $f_1, f_2$ do not coincide then we get the unique criterion 
$\arg\min_{q\in F}\{f_1(q)+f_2(q)\}$. 
With the same choice  as above we get the table 2 wich shows that  $\mathcal{N}(0,1)$  is the best proposal distribution.  

\begin{table}[!ht]
\begin{center}
  \begin{tabular}{|c||ccc|}
    \hline
    \scriptsize{Proposal} & \scriptsize{$f_1(q)$}   & \scriptsize{$f_2(q)$}  & \scriptsize{$f_1(q)+f_2(q)$} \\
    \hline
    \scriptsize{$\sigma^2=1$} & \scriptsize{0.0326} & \scriptsize{0.0102} & \scriptsize{0.0428}  \\ 
     \hline
    \scriptsize{$\sigma^2=9$ } & \scriptsize{0.0338}  & \scriptsize{0.0155} & \scriptsize{0.0493} \\ 
      \hline
     \scriptsize{$\sigma^2=16$ } & \scriptsize{0.0378}  & \scriptsize{0.0188} & \scriptsize{0.0559} \\ 
      \hline
\end{tabular}
  \caption{$y=t=1$, $T=0.1$,  $N=5000$,$M=600$.}
  \end{center}
\end{table}

3) The case $y >t$. 

a) Criterion using the asymptotic bias: We propose $X_T(y,t)$ as an estimator of $soft(y,t)=y-t$. The mean $\Eb[X_T(y,t)-(y-t)]\approx 0$ for small $T$. In order to take account of all $y >t$ we draw $y$
from $Pareto(\alpha, t)$ distribution. We showed that the best choice is $\alpha=3$. 
Let $(X_i: i=1, \ldots M)$
be a sample of $Pareto(\alpha, t)$ with $\alpha=3$. The best proposal distribution for sampling $X_T(y,t)$ will produce a sequence $(\theta^{(n)}(y): \quad n=1, \ldots, N)$ such that 
\be 
f_1(q)=\frac{1}{M}\sum_{i=1}^{M}|\frac{1}{N}\sum_{n=1}^{N}(\theta^{(n)}(y)-(y-t))|
\ee 
is minimal. 

b)Criterion using the mean square error: The mean square error $\Eb[(X_T(y,t)-(y-t))^2]$ for small $T$ is equal to 
$T\Eb[\mathcal{N}^2(0,t)]=Tt$. 
The best proposal distribution for sampling $X_T(y,t)$ will produce a sequence $(\theta^{(n)}:\quad n=1, \ldots, N)$ such that  
\be 
f_2(q)=\frac{1}{M}\sum_{i=1}^M|\frac{1}{N}\sum_{n=1}^N(\theta^{(n)}(X_i)-(X_i-t))^2-Tt|
\ee 
is minimal.
If the minimizers of $f_1$, $f_2$ do not coincide then we get the unique criterion
\be 
\arg\min_{q\in F}\{f_1(q)+f_2(q)\}.
\ee
According to Table 3, for $t=1$,  $\mathcal{N}(0,1)$ is also the best proposal distribution.  
\begin{table}[!ht]
\begin{center}
  \begin{tabular}{|c||ccc|}
    \hline
    \scriptsize{Proposal} & \scriptsize{$f_1(q)$} & \scriptsize{$f_2(q)$} & \scriptsize{$f_1(q)+f_2(q)$} \\
    \hline
    \scriptsize{$\sigma^2=1$} & \scriptsize{0.1388}& \scriptsize{0.0204} & \scriptsize{0.1591} \\ 
      \hline
    \scriptsize{$\sigma^2=9$} & \scriptsize{0.1397}& \scriptsize{0.0222} & \scriptsize{0.1619 } \\
      \hline
    \scriptsize{$\sigma^2=16$ } & \scriptsize{0.1419}& \scriptsize{0.0233} & \scriptsize{0.1652} \\
       \hline
\end{tabular}
  \caption{ $t=1$, $y>t$, $T=0.1$,  $N=5000$, $M=600$.}
  \end{center}
\end{table}\\

\section{Choice of the temperature in Metropolis-Hasting's algorithm}  
In this section we discuss the temperature needed in the estimation of    
the PLSE $soft(y,t)$ using our adaptative Metropolis Hasting's 
algorithm. The idea is to fix the bias $b$ and the mean square error $MSE$, and 
then choose the temperature $T$ such that $\Eb[X_T(y,t)]\approx b$, $\Eb[X_T^2(y,t)]\approx MSE$.   
We distinguish three cases. 

1) The case $y\in (0,t)$. 
 
a) Controlling the asymptotic bias: Fixing for small $T$ the bias 
\be 
\E[X_T(y,t)]\approx T\E[X(y,t)]=T\frac{2\frac{y}{t}}{(1-\frac{y^2}{t^2})}
=Tm_1(\frac{y}{t}):=b,  
\ee 
we get, for $y\neq 0$, the temperature $T(b,\frac{y}{t}):=\frac{b}{m_1(\frac{y}{t})}$. 
We plot in Figure 2 (a), for $b=0.001$, $u\in (0,1)\to T(b,u)$.

b) Controlling the asymptotic mean square error: Fixing for small $T$ the mean square error $\Eb[X_T^2(y,t)]$ 

\be 
\Eb[X_T^2(y,t)]\approx T^2\Eb[X^2(y,t)]
=T^2\{\frac{1-\frac{y}{t}}{(1+\frac{y}{t})^2}+\frac{1+\frac{y}{t}}{(1-\frac{y}{t})^2}\}=T^2m_2(\frac{y}{t}):=MSE,  
\ee 
we get the temperature $T(MSE,\frac{y}{t})=\sqrt{\frac{MSE}{m_2(\frac{y}{t})}}$.
We plot in Figure 2 (b), for $MSE=0.01$, $u\in (0,1)\to T(MSE,u)$.

 \begin{figure}[H]
  \centering
     \includegraphics[width=8cm]{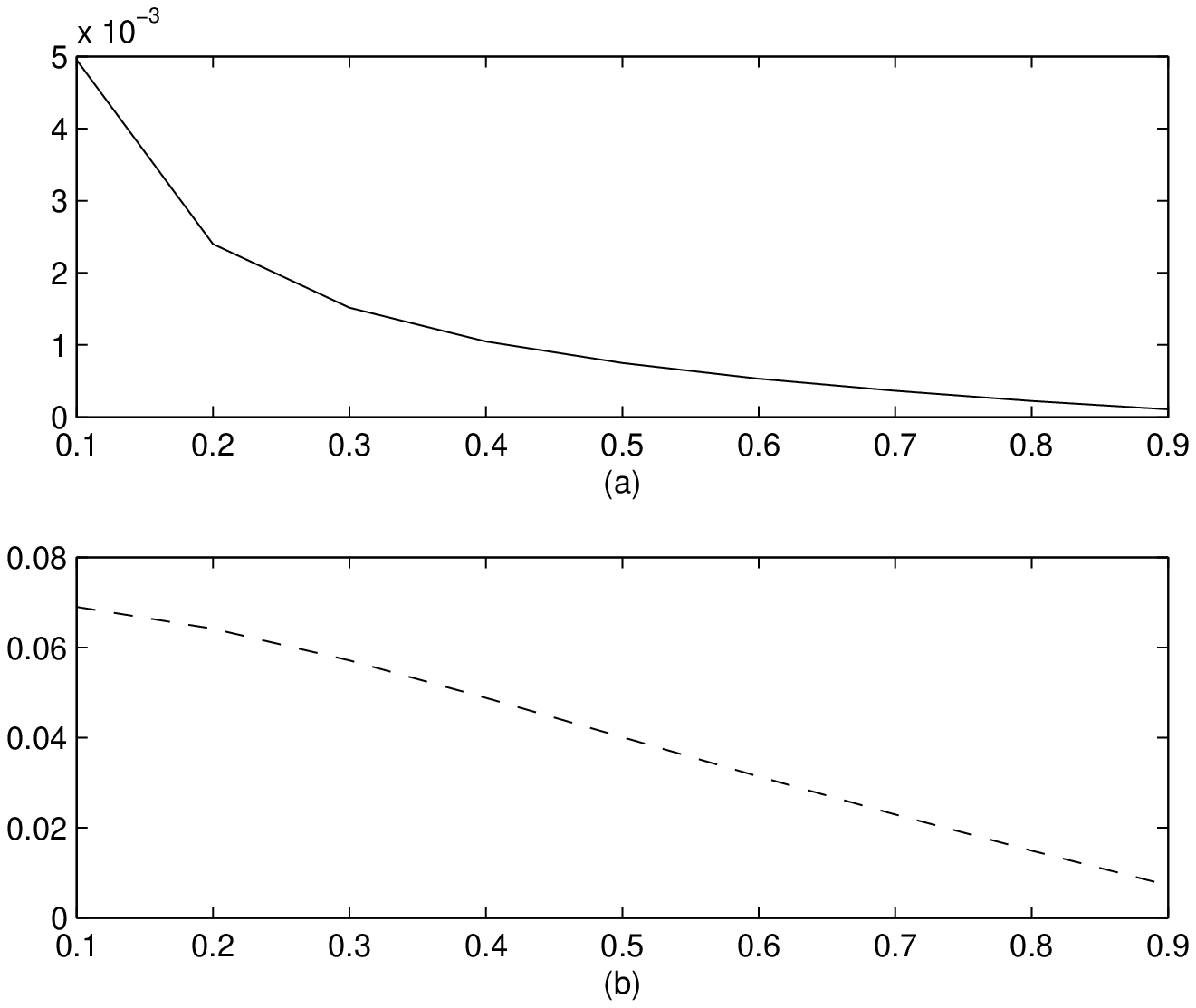}
 \caption{$b=0.001$, $MSE=0.01$, $u\in (0,1)\to T(b,u)$ and $u\in (0,1)\to T(MSE,u)$.}
\end{figure}

Now, we define our criterion for choosing the temperature as follows: $T_{b,MSE}(u):=T(b,u)=T(MSE,u)$.    
In order to have $T(b,\frac{y}{t})=T(MSE,\frac{y}{t})$, we need the constraint 
$\frac{b^2}{MSE}=\frac{m_1^2(\frac{y}{t})}{m_2(\frac{y}{t})}$ between the bias and the mean square error.
We plot in Figure 3 (a), (b), respectively the latter constraint as a function of $\frac{y}{t}\in (0,1)$ and 
the map $MSE\in (0,2)\to T(MSE,u)$ with $u=0.5$.   
 \begin{figure}[H]
  \centering
     \includegraphics[width=8cm]{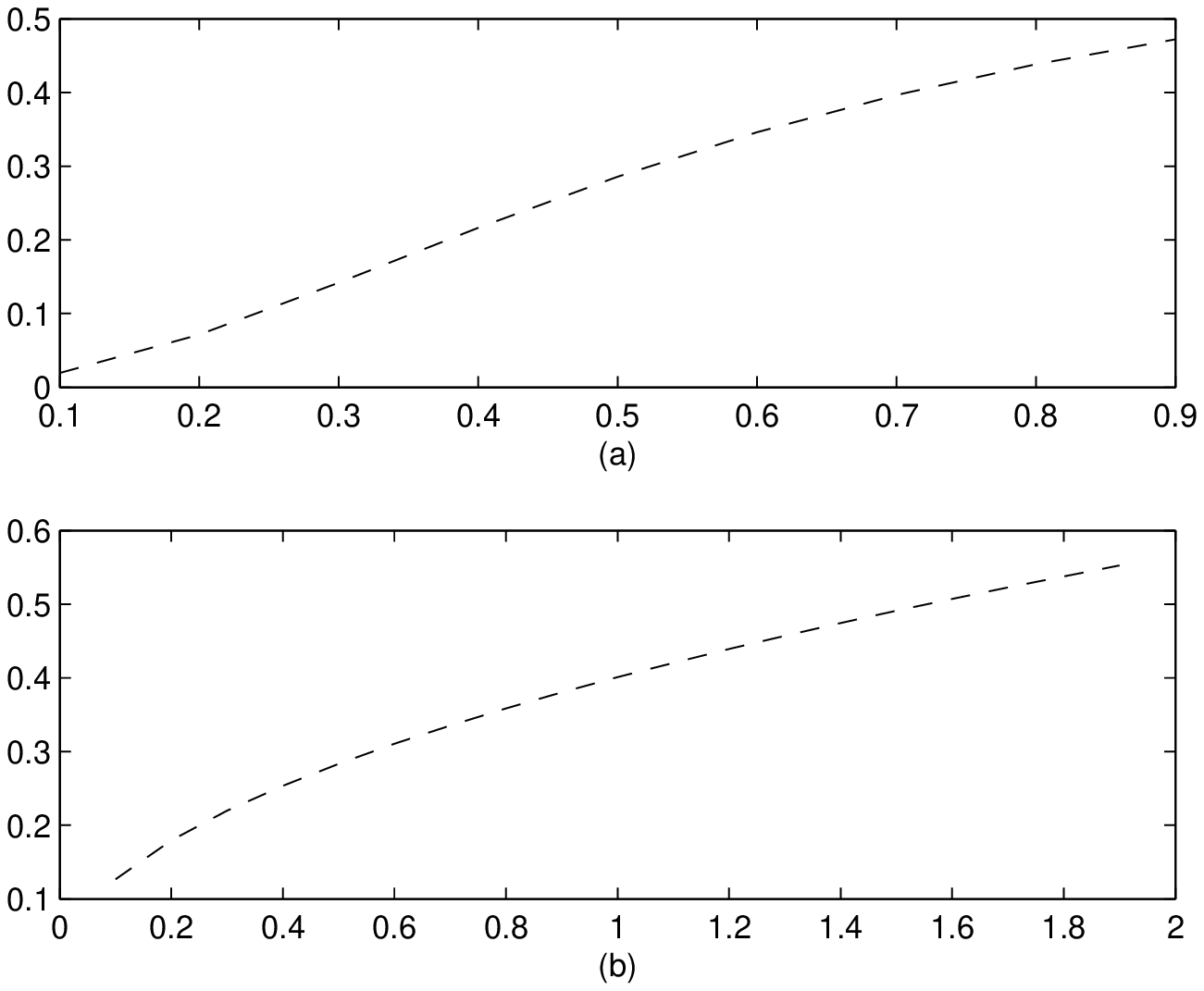}
 \caption{ $u=0.5$, $MSE\in (0,2)\to T(MSE,u)$ and  $u\in (0,1)\to \frac{m_1^2(u)}{m_2(u)}$ .}%
\end{figure}
\bigskip 

2) The case $y=t$. 

a) Controlling the asymptotic bias: Fixing for small $T$ the bias 
\be 
\Eb[X_T(t,t)]\approx\sqrt{T}\Eb[|\mathcal{N}(0,t)|]=\sqrt{\frac{2Tt}{\pi}}:=b,
\ee   
we get the temperature $T(b)=\frac{\pi b^2}{2t}$.

b) Controlling the asymptotic mean square error: Fixing for small $T$ the mean square error 
\be \Eb[X_T^2(t,t)]\approx T\Eb[N^2(0,t)]=Tt:=MSE,  
\ee 
 we get the temperature $T(MSE)=\frac{MSE}{t}$. 
In order to have the same temperature we set 
$T_{b,MSE}:=T(b)=T(MSE)$. The latter equality implies the relation $MSE=\frac{\pi b^2}{2}$ between the bias and the mean square 
error.  

3) The case $y >t$. Here the bias $b=0$ and we need only a Fixed mean square error, i.e.   
\be 
\Eb[(X_T(y,t)-(y-t))^2]\approx T\Eb[\mathcal{N}^2(0,t)]=Tt:= MSE,  
\ee 
we get the temperature $T_{0,MSE}=\frac{MSE}{t}$.

\subsection{Metropolis-Hasting's algorithm} 
Metropolis-Hasting's algorithm produces a Markov chain $(\theta_{MH}^n)$ such that for any suitable measurable function $h$  
\be 
\Eb[h(X_T(y,t)]=\lim_{N\to +\infty}\frac{\sum_{n=0}^Nh(\theta^n)}{N}. 
\ee  
In this section we address the problem of the convergence of the series 
$\frac{\sum_{n=0}^Nh(\theta^n)}{N}$
in the cases $h(x)=x$, $h(x)=x^2$ and $y\in (0,t)$. 
We fix the bias $b$ and the corresponding mean square error $MSE$. We derive 
the temperature $T_{b,MSE}$ solution of $\Eb[X_T(y,t)]=b$ and $\Eb[X_T^2(y,t)]=MSE$. We run Metropolis-Hasting's algorithm with the temperature $T_{b,MSE}$ and we calculate 
the sums $b_N:=\frac{1}{N}\sum_{n=0}^N\theta^n$, and $MSE_N:=\frac{1}{N}\sum_{n=0}^N(\theta^n)^2$  for different sizes $N$. 
In The table 4 we fix $b=0.01$, $MSE=0.00035$. We vary $N$ and we calculate the values of $b_N$ and $MSE_N$.
We show for $N=8000$, that $b_N\approx b$ and $MSE_N\approx MSE$.          
\begin{table}[!ht]
\begin{center}
  \begin{tabular}{|c||ccc|}
   \hline
   \scriptsize{$N_{MH}$} & \scriptsize{2000} &  \scriptsize{5000}  &  \scriptsize{8000} \\
    \hline
    \scriptsize{$b_N$} & \scriptsize{0.0066} & \scriptsize{0.0114} & \scriptsize{0.0106} \\ 
     \hline
    \scriptsize{$MSE_N$} & \scriptsize{6.2532e-04}  & \scriptsize{5.905e-04} & \scriptsize{3.541e-04} \\
      \hline
    
\end{tabular}
  \caption{  $b_N$ and $MSE_N$ values , for  $y=0.5$, $t=1$,  $b=0.01$, $MSE=3.5e-04$}
  \end{center}
\end{table} 
\subsection{Simulated annealing algorithm and comparison with Metropolis-Hasting's 
algorithm} 
We address the convergence of simulated annealing's algorithm to $soft(y,t)$ with $y\in (0,t)$. We consider the geometric tempering $\beta_n=\beta_{0}q^n=\frac{1}{T_n}$, with $q=1.001$ and $\beta_0=1$ . 
Fixing the bias $b$ and the corresponding mean square error $MSE$ we get the iteration number $N_{b,MSE}(SA)$ of simulated annealing algorithm to reach the temperature 
$T_{b,MSE}$. The number $N_{b,MSE}(SA)$ is the solution of the equation $\beta_{0}q^n=\frac{1}{T_{b,MSE}}$ i.e. $N_{b,MSE}(SA)=\frac{\ln(\frac{T_0}{T_{b,MSE}})}{\ln(q)}$.
Now, we compare the means of $\theta_{MH}^{N_{b,MSE}(MH)}$ and $\theta_{SA}^{N_{b,MSE}(SA)}$. Here $(\theta^n_{MH})$ and $(\theta^n_{SA})$ denote the sequences produced respectively 
by Metropolis-Hasting's and simulated annealing algorithms. We plot in Figure 4 (a), (b) the map $MSE\in (0,0.1)\to N_{b,MSE}(SA)$ respectively for $y=0$, $t=1$, and $y=0.5$, $t=1$.  
 We plot the maps $n=1...N_{b,MSE}(SA)\to \theta_{MH}^n$, $n=1...N_{b,MSE}(SA)\to \theta_{SA}^n$ respectively in Figure 5 (a), (b). 
\begin{figure}[!ht]
  \centering
     \includegraphics[width=8cm]{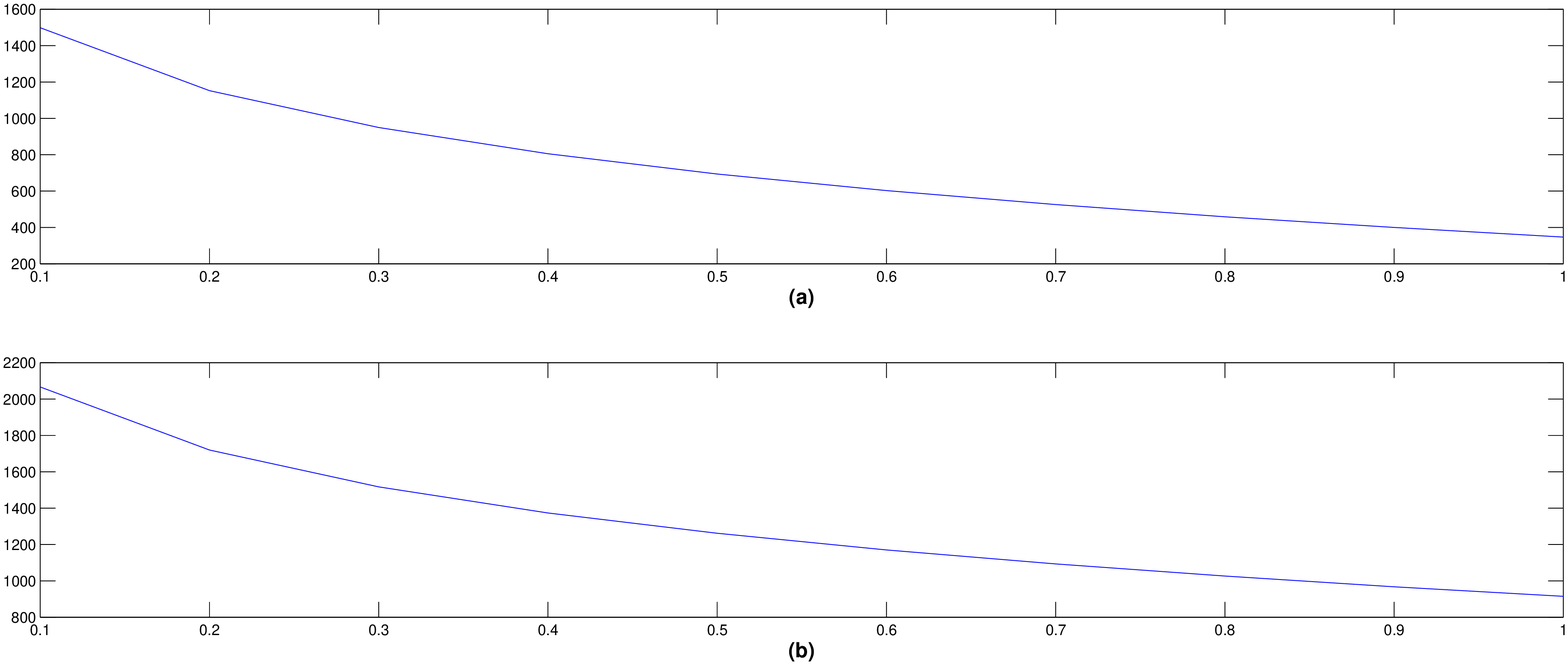}
 \caption{ $MSE\in (0,0.1)\to N_{b,MSE}(SA)$ .}
\end{figure}

\begin{figure}[!ht]
 \centering
    \includegraphics[width=8cm]{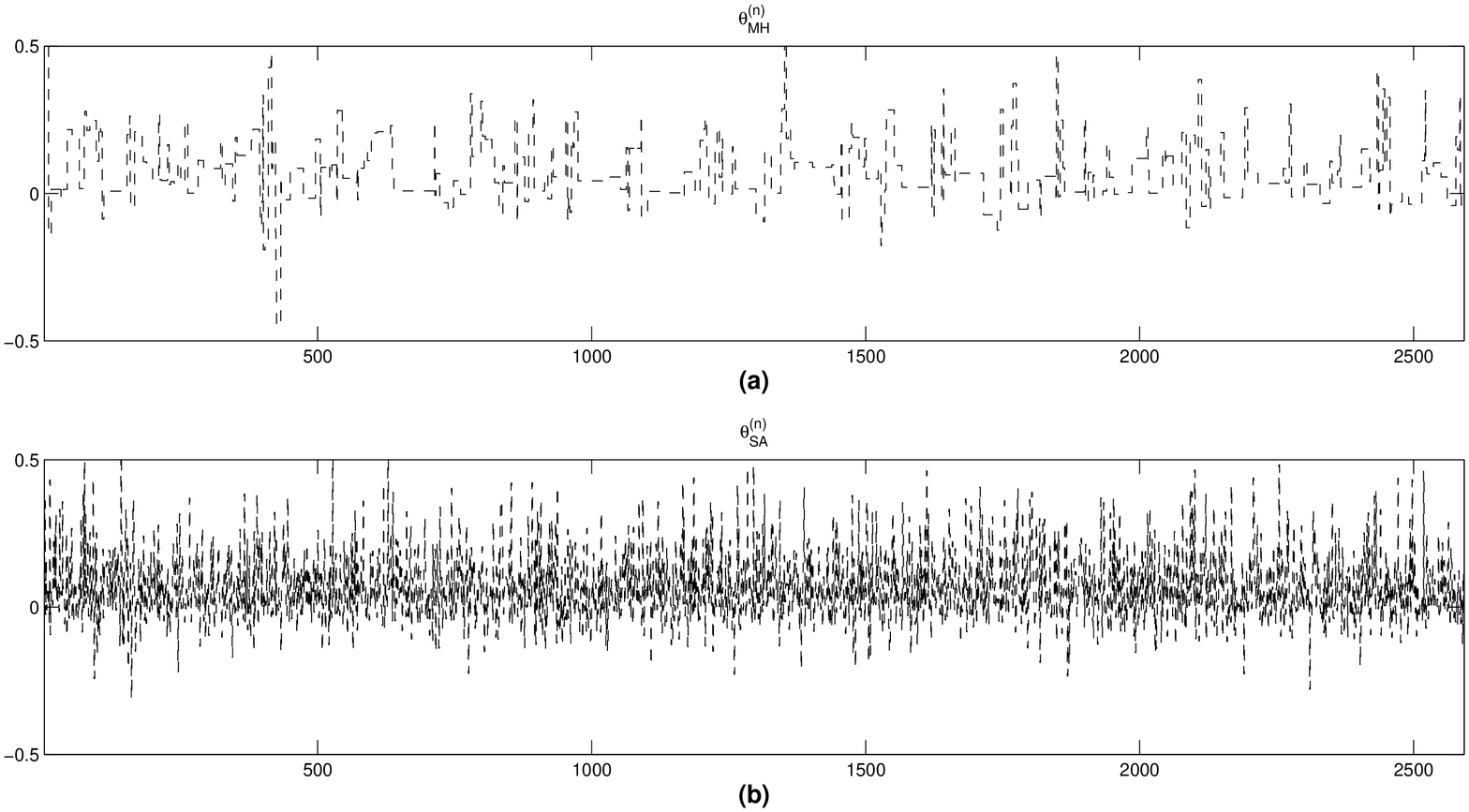}
 \caption{ $n\to \theta_{MH}^n$ and $n\to \theta_{SA}^n$.}
\end{figure}

\section{Conclusion}
In this paper we treated the Basis Pursuit De-noising problem
using Gibbs measures. We obtained the scaling of these Gibbs measures as the temperature goes to zero. We got, thanks to this scaling,    
several criteria to choose proposal distribution to initialize the Metropolis-Hasting's algorithm, and new criteria for choosing the temperature.
We also compared Metropolis-Hasting's and  simulated annealing algorithms. Our results can be easily extended to the  
analysis sparsity problem i.e. the minimization of the objective function 
$\|\D\x\|_1+\frac{\|\A\x-\y\|^2}{2t}$ with $\D\neq \I$.

\end{document}